\documentclass{IEEEtran}
\IEEEoverridecommandlockouts
\usepackage{cite}
\usepackage{amsmath,amssymb,amsfonts}
\usepackage{url}
\usepackage{algorithmic}
\usepackage{graphicx}
\usepackage{textcomp}
\usepackage{xcolor}
\def\BibTeX{{\rm B\kern-.05em{\sc i\kern-.025em b}\kern-.08em
    T\kern-.1667em\lower.7ex\hbox{E}\kern-.125emX}}

\usepackage{tikz}
\usepackage{xcolor}
\usetikzlibrary{patterns}
\usepackage{algorithm}
\usetikzlibrary{arrows.meta,positioning} 
\usepackage{thmtools}
\usepackage{thm-restate}
\usepackage{amssymb}
\usetikzlibrary{decorations.pathreplacing}

\usepackage{pgfplots}
\pgfplotsset{compat=1.12}
\usepgfplotslibrary{patchplots}
\usepgfplotslibrary{statistics}
\usepgfplotslibrary{external} 
\usepgfplotslibrary{fillbetween}
\usepackage{subfig}
\usepackage{caption}
\usetikzlibrary{external}
\usepackage{pgfplotstable}

\pgfplotsset{select coords between index/.style 2 args={
    x filter/.code={
        \ifnum\coordindex<#1\fi
        \ifnum\coordindex>#2\fi
    }
}}

\definecolor{colorpyC0}{RGB}{31, 119, 180} 
\definecolor{colorpyC1}{RGB}{255, 127, 14} 
\definecolor{colorpyC2}{RGB}{44, 160, 44} 
\definecolor{colorpyC3}{RGB}{214, 39, 40} 
\definecolor{colorpyC4}{RGB}{148, 103, 189} 
\definecolor{colorpyC5}{RGB}{140, 86, 75} 
\definecolor{colorpyC6}{RGB}{227, 119, 194} 
\definecolor{colorpyC7}{RGB}{127, 127, 127} 
\definecolor{colorpyC8}{RGB}{188, 189, 34} 
\definecolor{colorpyC9}{RGB}{23, 190, 207} 

\newcommand{\focs}{\textsc{Focs}}
\newcommand{\qp}{\textsc{QP}}
\newcommand{\focsmpc}{\focs-pm}
\newcommand{\mpc}{\textsc{Mpc}}

\newtheorem{thm}{Theorem}
\newtheorem{lem}[thm]{Lemma}

\begin{document}

\title{Improving the Optimization in Model Predictive Controllers: Scheduling Large Groups of Electric Vehicles}

\author{Leoni Winschermann, Marco E.T. Gerards, Johann Hurink\\
        Department of Electrical Engineering, Mathematics and Computer Science\\
        University of Twente\\
        Enschede, The Netherlands\\
        L.Winschermann@utwente.nl}

\date{}

\maketitle
\begin{abstract}
In parking lots with large groups of electric vehicles (EVs), charging has to happen in a coordinated manner, among others, due to the high load per vehicle and the limited capacity of the electricity grid. 
To achieve such coordination, model predictive control can be applied, thereby repeatedly solving an optimization problem. Due to its repetitive nature and its dependency on the time granularity, optimization has to be (computationally) efficient.

The work presented here focuses on that optimization subroutine, its computational efficiency and how to speed up the optimization for large groups of EVs. 
In particular, we adapt \focs, an algorithm that can solve the underlying optimization problem, to better suit the repetitive set-up of model predictive control by adding a pre-mature stop feature. 
Based on real-world data, we empirically show that the added feature speeds up the median computation time for 1-minute granularity by up to 44\%. 
Furthermore, since \focs\ is an algorithm that uses maximum flow methods as a subroutine, the impact of choosing various maximum flow methods on the runtime is investigated.
Finally, we compare \focs\ to a commercially available solver, concluding that \focs\ outperforms the state-of-the-art when making a full-day schedule for large groups of EVs.
\end{abstract}


\section{Introduction}
With the on-going transition to electric mobility, the number of electric vehicles (EVs) in the Netherlands is increasing rapidly~\cite{2020WolbertusPHEVadoptionNLlessonsLearned}. 
Especially near-work locations have been found to be good locations to charge large groups of EVs, partially due to the long stays~\cite{2018SadeghianpourhamamiQuantitativeAnalysisofEVflexibilityDataDrivenApproach} and the potential availability of solar energy~\cite{2021OsorioRooftopPVParkingLotsSupportEVChargingReview}. 
However, due to limited grid capacity~\cite{2014EisingSmartGridsRisksIntegrationEnergyandTransportation}, high synchronicity~\cite{2010TuritsynRobustBroadcastCommControlofEVCharging} and possible temporal mismatches with said solar energy production~\cite{2021CanigueralEVprofileClusteringArnhemCaseStudy}, charging has to occur in a coordinated manner.
For successful coordination, a controller needs information on the charging sessions, for example the (expected) energy demand and departure time of an EV in the parking lot.
In practice, the EVs currently do not communicate such information~\cite{2019VecchioMayTheForceMoveYouRolesActorsofInformationSharingDevices}.
Therefore, next to the sheer amount of charge required per vehicle, the biggest challenge in managing the charging processes of the EVs is the absence of information.

To bridge this information gap, model predictive control (\mpc) can be applied (e.g.,~\cite{2021ZhouQuantitativeMPCasEMS4hybridEV}). Based on the information available at the current point in time and an underlying system model, \mpc\ generates and solves a suitable offline optimization problem and an (in time) initial portion of this solution is realized. This process is repeated frequently over the given time horizon. 
Specifically to solve this offline problem, the Flow-Based Offline Charging Scheduler (\focs) has been developed~\cite{2023WinschermannFOCSArxiv}. \focs\ is a polynomial time algorithm which in theory has a complexity of $\mathcal{O}(n^2\mu)$, where $n$ is the number of EVs to be charged and $\mu$ is the complexity of the maximum flow algorithm embedded in \focs. 

Notably, research has shown that for some algorithms theoretical and empirical runtime do not align. For example, the simplex algorithm, while theoretically exponential in runtime, performs well in practice. The other way around, Ellipsoid methods for linear programs are known to be polynomial on paper, while performing poorly in practice. 
This motivates us to investigate the empirical runtime of the novel algorithm \focs\ based on real world EV data.

The main contributions of the work presented here are:
\begin{itemize}
    \item Introduction and validation of a pre-mature stop feature that increases computational efficiency in model predictive controllers using \focs.
    \item Extensive comparative runtime analysis of \focs, and its validation using real world EV charging data. 
\end{itemize}

The remainder of the paper is organized as follows. Section~\ref{sec:EVproblembackground} discusses the mathematical model of the considered EV scheduling problem (Section~\ref{sec:problemdefintion}), briefly introduces \focs\ (Section~\ref{sec:focsNutshell}) and proposes a pre-mature stop feature for \focs\ (Section~\ref{sec:introfocsmpc}) expected to increase its efficiency for \mpc\ applications. Then, Section~\ref{sec:methods} specifies the experimental setup and gives some information on the implementation of the used algorithms. Furthermore, the section describes the real world data set used for the experiments (Section~\ref{sec:datasetACM}). Section~\ref{sec:results} presents the experimental results, followed by a discussion in Section~\ref{sec:discussion} and a conclusion in Section~\ref{sec:conclusion}.

\section{EV scheduling problem and algorithms} \label{sec:EVproblembackground}
In this section, we describe the considered EV scheduling problem in terms of notation, constraints and considered objective function. Furthermore, we briefly introduce \focs, an algorithm that can be used to determine an optimal solution for said problem. For more details on \focs, we refer to~\cite{2023WinschermannFOCSArxiv}.

\subsection{Problem definition} \label{sec:problemdefintion}
We consider jobs $1, ..., n$, where each job $j$ corresponds to a pending charging session (or EV) with associated energy requirement $e_{j}$, arrival time $a_j$, departure time $d_j$ and job-specific charging power limit $p_j^{\max}$. We denote the set of all jobs by $[n]$.
We consider a problem instance feasible if 
\begin{displaymath}
    e_j \leq p_j^{\max}(d_j - a_j) \hspace{10pt} \forall j \in [n].
\end{displaymath}

Moreover, we discretize the time horizon $[ \min_j a_j, \max_j d_j]$ into atomic intervals by using all arrival and departure times of EVs as breakpoints. The resulting (ordered) sequence of breakpoints we denote as $t_1, \dots, t_{m+1}$, resulting in $m \leq 2n -1$ atomic intervals $I_i = [t_i, t_{i+1}]$  for $i = 1, ..., m$.

To relate jobs to atomic intervals, we denote by $J(i)$ the jobs that are available in interval $I_i$, i.e., 
\begin{displaymath}
        J(i) =\{ j | (a_j \leq t_i) \land (t_{i+1} \leq d_j)\}.
\end{displaymath}
Similarly, $J^{-1}(j)$ is defined as the set of indices $i$ of the intervals $I_i$ where $j$ is available.
Finally, as atomic intervals are in general not unit-sized, we introduce maximum energy limits $e_{i,j}^{\max} = p_j^{\max} (t_{i+1} - t_i)$ per job $j$ and interval $I_i$.

As decision variables, let $e_{i,j}$ be the energy that EV $j$ charges in $I_i$.
Note that preemption is allowed, meaning the charging can be suspended for a few intervals to be continued at a later interval. In this work, we do not consider V2G, implying that $e_{i,j}\geq 0$. 

Summarizing, the EV model is constrained by
\begin{subequations}\label{eq:MIP}
\begin{align}
    \sum_{i\in J^{-1}(j)} e_{i,j} &\geq e_{j}  &\forall j\in[n] \label{eq:MIPnoENS}\\
    e_{i,j} &\geq 0& \forall j \in [n], i\in J^{-1}(j) \hspace{3pt}\label{eq:MIPnonnegLoad}\\
    e_{i,j} &\leq e_{i,j}^{\max} & \forall j \in [n], i\in J^{-1}(j). \label{eq:MIPspeedlimit}
\end{align}
\end{subequations}

From a grid perspective, the aggregated power level resulting from an EV schedule is of interest. For a given schedule, the average aggregated power level $p_i$ in atomic interval $I_i$ is given by
\begin{displaymath}
        p_i = \frac{\sum_{j\in J(i)} e_{i,j}}{t_{i+1} - t_i}. 
\end{displaymath}

The two most frequently considered objectives for the overall grid usage are the $\ell_{\infty}$ and $\ell_{2}$-norms of the aggregated power. \focs\ was developed for minimization of objective functions 
\begin{equation}
        F\left( p_1, \dots, p_m\right) \label{eq:MIPobjective}    
\end{equation}
that are convex, differentiable and for which increasing any value $p_i$ increases the value of the objective function. The complexity of \focs\ is independent on the exact form of the objective function. In this work, we specifically consider the square of the $\ell_2$-norm given by
\begin{equation}
    F\left( p_1, \dots, p_m\right) = \sum_{i = 1}^m L_i \cdot \left(p_i^2\right) \label{eq:L2NormObjective}
\end{equation}
where $L_i$ is a normalization term to account for the length of the interval. 

\subsection{The general working of \focs} \label{sec:focsNutshell}
\focs\ is a recursive algorithm that finds an optimal solution to the problem described above. To this end, it makes use of various properties of optimal solutions to this problem. First, any optimal aggregated power profile that minimizes an objective function (\ref{eq:MIPobjective}) is a step function in the sense that within each atomic interval, the aggregated power is constant. Based on this, \focs\ iteratively identifies a critical interval where the highest speed is needed, determines the minimum speed for this interval, and fixes the schedule for this interval, ultimately minimizing the maximum peak. The selected critical interval and the charge provided in that interval are then removed from the problem instance. 

\begin{figure}[]
    \centering
    \includegraphics{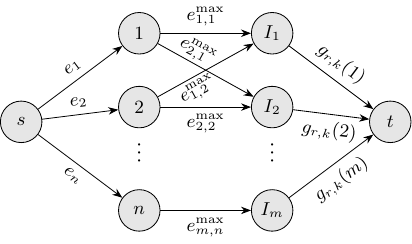}
    \caption{Schematic of flow network structure of EV charging schedule.}

    \label{fig:flowNetwork}
\end{figure}

To determine the critical intervals, \focs\ computes multiple iterations of maximum flows. Figure~\ref{fig:flowNetwork} depicts the basic structure of the considered flow network. Here, at each iteration, \focs\ adapts the capacities of edges adjacent to sink node $t$. The iterations until a critical interval is identified, form a round $r$. The edge capacities $g_{r,k}:\{I_i|i\in[m]\}\rightarrow\mathbb{R}_{\geq0}$ in round $r$ and iteration $k$ represent a fill-level that is increased over the iterations until critical interval(s) are found.

\subsection{Adapting \focs\ for model predictive control} \label{sec:introfocsmpc}
\focs\ was developed to determine the optimal solution of the offline EV scheduling problem described in Section~\ref{sec:problemdefintion}. To this end, it solves the instance over the entire planning horizon. However, for \mpc, often only the next (few) time intervals and their control actions are of interest. Based on this, we propose an adaptation to \focs. For convenience, we assume in the following that only the first time interval is of interest. The proposed concept extends naturally to include the first few intervals.

\begin{figure}[]
    \centering
    \begin{tabular}{c|c c c c}
         $j$ & $a_j$ & $d_j$ & $e_j$ & $p_j^{\max }$ \\  \hline
         1 & 0 & 2 & 4 & 2 \\
         2 & 1 & 3 & 2 & 1 \\
    \end{tabular}
    \\ \vspace{6 pt}
	\includegraphics{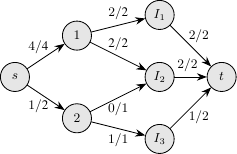}	
	\includegraphics{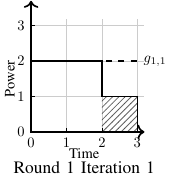}\\
	
	\includegraphics{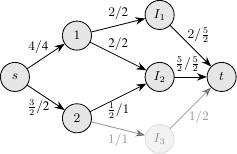}	
	\includegraphics{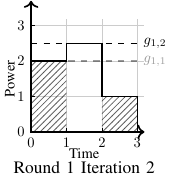}\\
	
	\includegraphics{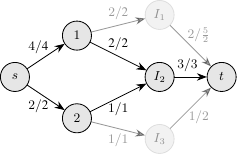}	
	\includegraphics{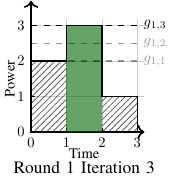}\\

	\includegraphics{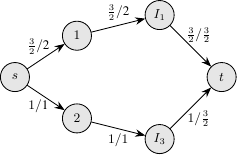}	
	\includegraphics{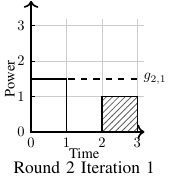}\\

	\includegraphics{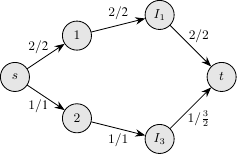}	
	\includegraphics{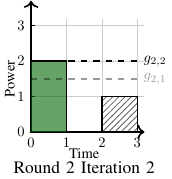}\\

	\includegraphics{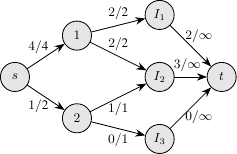}	
	\includegraphics{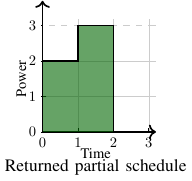}\\

    \vspace{6pt}
    \caption{Intermediate states of \focs \ with pre-mature stop for an example instance, tracked over rounds and iterations.}
    \label{fig:flowExWithProfileFOCSmpc}
\end{figure}

\focs\ recursively identifies intervals that require the highest aggregated power in an optimal schedule (so-called critical intervals), determines the partial schedule corresponding to that interval, and repeats this with the remaining problem instance. In particular, each of those recursion steps fixes a part of the optimal schedule. Therefore, if we are only interested in the schedule of the first interval, we can pre-maturely terminate \focs\ once the schedule for that interval was found, without losing any information on the first interval. 

We illustrate this concept in Figure~\ref{fig:flowExWithProfileFOCSmpc} for an instance with two EVs and three intervals. The annotations at the edges of the flow network are of the form $x/y$, where $x$ is the flow going through the edge, and $y$ is the edge capacity. 
Moreover, the green shaded blocks are the identified critical intervals. 
Those are fixed, and in the last step combined into the algorithm output. After finding the schedule for the first interval (i.e., when we find the green block for the first interval in the fifth step in Figure~\ref{fig:flowExWithProfileFOCSmpc}), we terminate the algorithm. Note that even though the schedule for the third interval has not yet been found and thereby the returned schedule (last power profile in Figure~\ref{fig:flowExWithProfileFOCSmpc}) is not feasible considering the entire problem instance, there exists an optimal schedule that complements the found partial solution. In particular, the aggregated power in the first interval is optimal. This follows by optimality of \focs, and the following lemmas paraphrased from~\cite{2023WinschermannFOCSArxiv}. 

\begin{lem} (Lemma 2 in~\cite{2023WinschermannFOCSArxiv})
    If a critical interval $I_i$ and its aggregated power are found, and if there is at least one interval $I_{i'}$ left that at that point has not yet been scheduled, then no scheduled charging work can feasibly be moved from $I_i$ to $I_{i'}$ for any such $I_{i'}$.
\end{lem}

\begin{lem} (Lemma 3 in~\cite{2023WinschermannFOCSArxiv})
    The aggregated power of the subsequently identified (sets of) critical intervals by \focs\ is strictly decreasing. In particular, if \focs\ fixes the schedule for $I_i$ before that of $I_{i'}$, then 
    \begin{displaymath}
                p_i > p_{i'} \nonumber
    \end{displaymath}
    in the optimal solution.
\end{lem}

For the precise formulation and proof of the lemmas, we refer the reader to \cite{2023WinschermannFOCSArxiv}. 

Notably, for EV parking lots near office locations, the optimal schedule over a day typically resembles a concave curve. The pre-mature stop feature we propose here terminates the algorithm once the first interval has been scheduled, where intervals are scheduled by order of their aggregated power in the optimal solution. Therefore, the time of day at which we apply the pre-mature stop feature impacts its effectiveness. In particular for an office parking lot, we expect the greatest reduction in computation time with pre-mature stop compared to the unrestricted \focs\ to be around noon.
In the following, we refer to the implementation of \focs\ with pre-mature stop as \focsmpc.

\section{Methods} \label{sec:methods}
In this section, we provide some details on the implementation of the three tested optimization models (Quadratic program in Gurobi (\qp), \focs\ and \focs\ with pre-mature stop (\focsmpc)), after which we describe the experimental setup and real-world data set used.\footnote{The source code with all implementations and the experimental setup can be found under \url{https://github.com/lwinschermann/FlowbasedOfflineChargingScheduler.}}

\subsection{\qp\ implementation}\label{sec:methodsQPimplementation}
Based on the problem constraints (\ref{eq:MIP}) and the quadratic objective function (\ref{eq:L2NormObjective}), we use the commercial solver Gurobi~\cite{gurobi} to solve this problem. Gurobi is used for comparison with state-of-the art solvers. Note that while \focs\ leads to optimal solutions also for more general objective functions, for practical reasons we compare its runtime to a Gurobi model and for that we restrict us to a quadratic objective function. 

\subsection{\focs\ implementation}\label{sec:methodsFOCSimplementation}
The \focs\ implementation used for the empirical runtime analysis of \focs\ has been developed specifically as a proof of concept for this paper. The code is written in python, heavily relying on the \texttt{networkx} package~\cite{2008Networkx}. This package is chosen for its user friendly flow models, and because various maximum flow algorithms are readily available within the package. An overview of the considered algorithms is provided in Table~\ref{tab:maxflowComplexity}.
\begin{table*}[]
    \centering
    \caption{Complexity table for four maximum flow algorithms. Here, $N$, $M$ and $U$ are respectively the number of nodes, edges and the maximum capacity in a given flow network.}
\begin{tabular}{l|c|c|c}
     Max flow algorithm & Complexity~\cite{2023CruzMejiaSurveyExactAlgsMaxFlowAndMinCostFlow} & Networkx complexity & Complexity on \focs\ network \\\hline 
     \texttt{shortest\_augmenting\_path()} & $\mathcal{O}(N M U)$ & $\mathcal{O}(N^2 M)$ & $\mathcal{O}(n^4)$ \\ 
     \texttt{edmonds\_karp()} & $\mathcal{O}(M^2N)$ or $\mathcal{O}(m^2\log U)$ & $\mathcal{O}(M^2N)$ & $\mathcal{O}(n^5)$ \\ 
     \texttt{preflow\_push()} & $\mathcal{O}(MN\log (\frac{N^2}{M}))$ & $\mathcal{O}(N^2\sqrt{M})$ & $\mathcal{O}(n^3)$ \\ 
     \texttt{dinitz()} & $\mathcal{O}(N^2M)$ or $\mathcal{O}(MN\log U)$ & $\mathcal{O}(N^2M)$ & $\mathcal{O}(n^4)$ \\ 
\end{tabular}         
\label{tab:maxflowComplexity}
\end{table*}
Next to that, we also implemented the adapted \focs\ with pre-mature stop as introduced in Section~\ref{sec:introfocsmpc} by adding a feature that stops the algorithm once the relevant parts of the schedule are computed. 
In the following, we refer to \focsmpc\ for results achieved with this pre-mature stop.

\subsection{Dataset}\label{sec:datasetACM}
The experiments described in this paper use data collected at a real life EV charging parking lot in Utrecht, the Netherlands~\cite{Nijenhuis2022}. The data was collected between September 1st 2022 and August 31st 2023, resulting in a total of 13694 charging sessions. Each session specifies (among others) an arrival time, departure time and total energy charged during the session. At the parking lot, most charging stations have two sockets each. If only one socket is occupied, the maximum charging power is 22~kW. Else, it is 11~kW per socket. Therefore, for the experiments, we set the maximum charging power per session to either 11~kW or 22~kW, depending on whether the average observed charging power is below or above 11~kW. From the 13694 sessions in the data, a maximum of 113 sessions was recorded in a single day.

\subsection{Numerical experiments}
The main focus of this paper is on optimization in \mpc\ for EV charging parking lots and in particular its efficiency. To this end, we present numerical experiments based on real world data, evaluating the empirical runtime of \focs\ and \focsmpc\ compared to the benchmark results by Gurobi. 

The main parameters to classify the experiments are the
\begin{itemize}
    \item instance size $n$, i.e., the number of EV charging sessions;
    \item time granularity, i.e., either 1 minute, 15 minutes, 30 minutes or 1 hour;
    \item used maximum flow method in \focs\ (and \focsmpc).
\end{itemize}
Given a set of parameters for the experiment, we randomly sample $n$ sessions from the dataset described in Section~\ref{sec:datasetACM} and apply all three optimization models. We record the CPU runtime using the function \texttt{time.process\char`_time()} in-built in the python package \texttt{time}. In particular, we record the CPU runtime for building the models (\qp, \focs, \focsmpc) and solving them.
For a given model and experiment, the total runtime is the sum of CPU time taken to build and solve the model.

As pointed out in Section~\ref{sec:introfocsmpc}, \focsmpc\ is most effective when the first interval is expected to have one of the highest power peaks in the optimal profile. Therefore, next to solving the entire instance over one day, we also repeat the experiment assuming optimization starting at noon. In that case, due to the characteristics of an office parking lot, the first interval is expected to have quite a high power. 
To get meaningful results, we repeat this process 500 times, starting at the instance sampling. The values reported in this work are the median runtimes over those 500 runs.
All experiments are run on an Intel Xeon E5-2630 v3 processor.

\section{Results} \label{sec:results}
In this section, we present and analyze the results of the runtime experiments. In particular, we are interested in the efficiency and therefore usability of the various algorithms for \mpc\ in (large) EV parking lots. To this end, we focus on the dependency of runtime on the instance size $n$. 
In Section~\ref{sec:resultsfocsmpc} we first discuss the pre-mature stop feature for \focs\ introduced in Section~\ref{sec:introfocsmpc}, and how it reflects in the empirical results. Then, in Section~\ref{sec:resultsMaxFlowComparison}, we focus on the impact that the different maximum flow methods have on the performance of \focs. Lastly, we compare the impact of various time granularities under invariant maximum flow methods in Section~\ref{sec:resultsRTcomparison}. 

\subsection{\focsmpc} \label{sec:resultsfocsmpc}
As discussed in Section~\ref{sec:introfocsmpc}, an \mpc\ approach is mostly interested in the next (few) control action(s). Therefore, we propose \focsmpc\ as an adapted \focs\ that speeds up computation by stopping the optimization pre-maturely when the control action for the first upcoming intervals is known. While by design we expect the feature to speed up computation, this section investigates the extend of the effect, and validates it empirically. 

\begin{figure*}[t]
\begin{center}
    \includegraphics{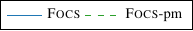}
\end{center}
\centering
    \subfloat[\centering 60s time granularity]{
     	\includegraphics{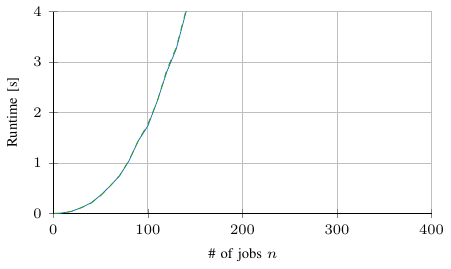}
    \label{fig:0_60_full_pmtest}
    }
    \subfloat[\centering 15m time granularity]{
    	\includegraphics{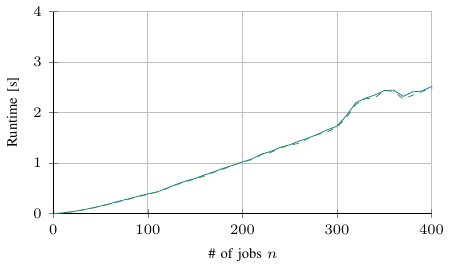}
    \label{fig:0_900_full_pmtest}
    } \\
    \subfloat[\centering 30m time granularity]{
		\includegraphics{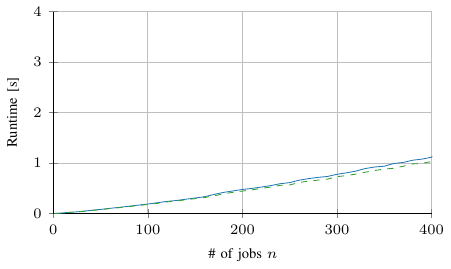}  
    \label{fig:0_1800_full_pmtest}
    }
    \subfloat[\centering 1h time granularity]{
		\includegraphics{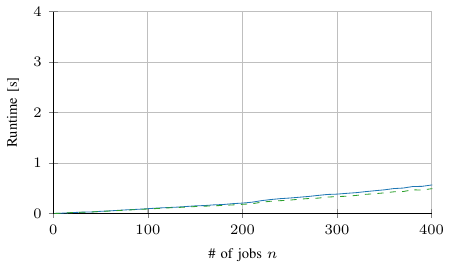}
    \label{fig:0_3600_full_pmtest}
    }
\caption{Runtimes relative to instance size, for various time granularities, solved for the whole day and using \texttt{shortest\_augmenting\_path()}.}
\label{fig:0_full_pmtest}
\end{figure*}

Figure~\ref{fig:0_full_pmtest} shows the total runtimes of \focs\ and \focsmpc\ relative to instance size $n$ for the four time granularities, and maximum flow method \texttt{shortest\_augmenting\_path()}. All instances are solved for the full day.
First, we note that the runtimes of \focs\ and \focsmpc\ are very similar for full day instances. Comparing the values over all instance sizes, the improvement of \focsmpc\ relative to \focs\ is at most 8\%, 11\%, 12\% and 15\% for time granularities of 1~minute, 15~minutes, 30~minutes and 1~hour respectively. The average improvements are 0\%, 2\%, 6\% and 9\% respectively. We see that especially for finer time granularities, the computational gain of the feature is rather small for this use case. 

In the analysis above, we considered instances for a full day. 
In a next step, we present results of experiments that consider instances starting at noon. 
\begin{figure*}[t]
\begin{center}
    \includegraphics{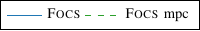}
\end{center}
\centering
    \subfloat[\centering 60s time granularity]{
        \includegraphics{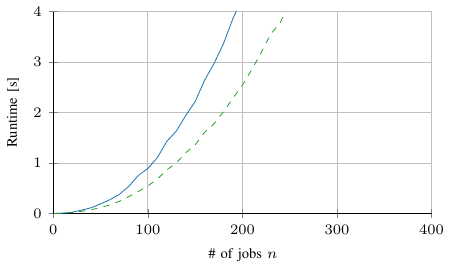}
    \label{fig:0_60_partial}
    }
    \subfloat[\centering 15m time granularity]{
    	\includegraphics{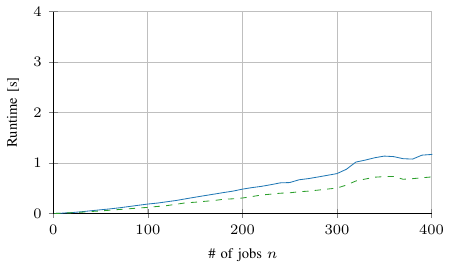}
    \label{fig:0_900_partial}
    } \\
    \subfloat[\centering 30m time granularity]{
		\includegraphics{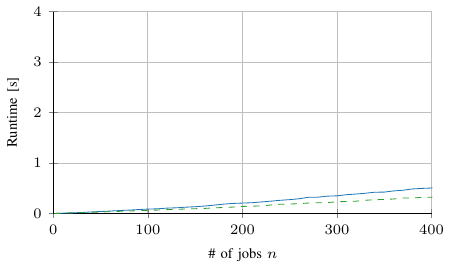}
    } \label{fig:0_1800_partial}
    \subfloat[\centering 1h time granularity]{
		\includegraphics{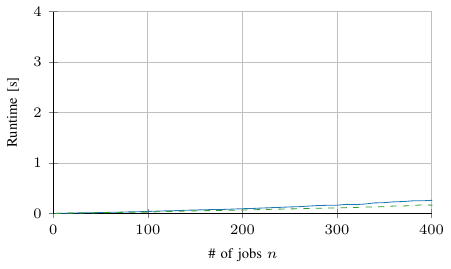}
    \label{fig:0_3600_partial}
    }
\caption{Runtimes for instances solved at noon using \texttt{shortest\_augmenting\_path()}.}
\label{fig:0_partial}
\end{figure*}
Figure~\ref{fig:0_partial} shows the total runtimes of \focs\ and \focsmpc\ relative to instance size $n$ for the four time granularities, and maximum flow method \texttt{shortest\_augmenting\_path()}. However, as opposed to Figure~\ref{fig:0_full_pmtest}, instances are solved starting from noon. 
As expected, for situations with the same instance size and time granularity, solving the (smaller) problem at noon is faster than solving the entire day. 
Furthermore, the added value of the pre-mature stop feature is even more pronounced in the noon experiments. For time granularities of 1~minute, 15~minutes, 30~minutes and 1~hour, the experiments show runtime improvements of at most 44\%, 39\%, 37\% and 38\% respectively for up to 400 EVs. That maximum improvement of the median runtime corresponds to instance sizes 400, 390, 290 and 340, all rather large instances.
The average improvements amount to 33\%, 28\%, 26\% and 23\% respectively.

\subsection{Maximum flow method comparison} \label{sec:resultsMaxFlowComparison}
The efficiency of an algorithm in the context of an \mpc\ among others depends on the implementation. As mentioned in Section~\ref{sec:focsNutshell}, \focs\ relies on solving multiple maximum flow problems.
Hereby, in theory it does not matter which concrete algorithm is used. However, since according to Table~\ref{tab:maxflowComplexity} the available maximum flow algorithms are non-linear, this section investigates their impact on the scalability of \focs\ for large groups of EVs. 

As mentioned in Section~\ref{sec:methodsFOCSimplementation}, \texttt{networkx} provides a number of functions that compute maximum flows. In particular, our experiments use the \texttt{shortest\_augmenting\_path()}, \texttt{edmonds\_karp()}, \texttt{preflow\_push()} and \texttt{dinitz()} functions from the \texttt{networkx} package. Table~\ref{tab:maxflowComplexity} provides an overview of those functions, and the following three properties. First, per function, we cite the theoretical complexity of its namesake based on \cite{2023CruzMejiaSurveyExactAlgsMaxFlowAndMinCostFlow}. Note that the implementations of \texttt{shortest\_augmenting\_path()}, \texttt{edmonds\_karp()}, \texttt{preflow\_push()} and \texttt{dinitz()} likely do not strictly follow the original algorithms (\cite{1956FordFulkersonMaxFlows}, \cite{1972EdmondsKarpMaxFlowAlgs}, \cite{1974KarzanovMaxFlowsPreflows}/\cite{1988GoldbergMaxFlows} and \cite{1970DinitzMaxFlows} respectively). Since their publication, various improvements and variants have been introduced that are usually referred to by the same name. Secondly, Table~\ref{tab:maxflowComplexity} presents the theoretical complexity as reported by \texttt{networkx}\footnote{last accessed 19 December 2023}. Finally, the last column takes the complexity reported by \texttt{networkx} and bounds the theoretical complexity specifically for flow networks of \focs\ instances.  
We note that the number of nodes (compare flow network in Figure~\ref{fig:flowNetwork}) is bound by 2 plus the number of jobs and intervals, being $2 + n + (2n-1) = 3n +1$. The number of edges is roughly bound by $n + n(2n-1) + (2n-1) = 2n^2 + 2n - 1$. We derive the theoretical complexity in terms of the input size of \focs\ instances by substituting this in the bound reported by \texttt{networkx}.
Since these maximum flow solvers are only used in the solving of \focs\ (and therefore do not influence the model building), we restrict the discussion of the results to the solving part of the runtime. 
    
\begin{figure}
	\includegraphics{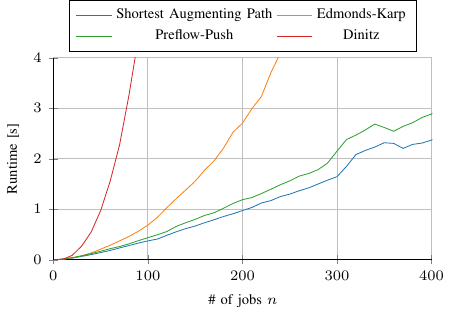}
  \caption{Runtimes for solving only for quarterly granularity, comparing maximum flow algorithms applied in \focs.}
  \label{fig:0123_900_full_sol}
\end{figure}

Figure~\ref{fig:0123_900_full_sol} presents the solving time of \focs\ for a full day and quarterly time granularity. 
Among the methods using augmenting paths, we see a clear difference in performance between methods. While the time taken to solve \focs\ using \texttt{shortest\_augmenting\_path()} seems to grow almost linearly for instance sizes up to 400, solving times with \texttt{edmonds\_karp()} and \texttt{dinitz()} clearly increase non-linearly. 
Finally, the preflow-push method behaves similarly to \texttt{shortest\_augmenting\_path()}.

If we compare this to the theoretical complexity of the various maximum flow methods embedded in \texttt{networkx} (see Table~\ref{tab:maxflowComplexity}), the most striking observation is that while \texttt{preflow\_push()} has the smallest theoretical runtime on \focs\ networks ($\mathcal{O}(n^3)$), the theoretically quartic runtime of \texttt{shortest\_augmenting\_path()} ($\mathcal{O}(n^4)$) outperforms the other tested methods. \texttt{Edmonds\_karp()} shows the highest theoretical complexity, and empirically is the second slowest method, overtaken only by \texttt{dinitz()}.

Overall, \texttt{shortest\_augmenting\_path()} appears to be the dominating method in our experiments with respect to runtime. The experiments further indicate that the proper choice of maximum flow method greatly increases the usability of \focs\ for EV scheduling.


\subsection{Runtime comparison}\label{sec:resultsRTcomparison}
\begin{figure*}[t]
\begin{center}
    \includegraphics{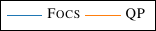}
\end{center}
\centering
    \subfloat[\centering 60s time granularity]{
     	\includegraphics{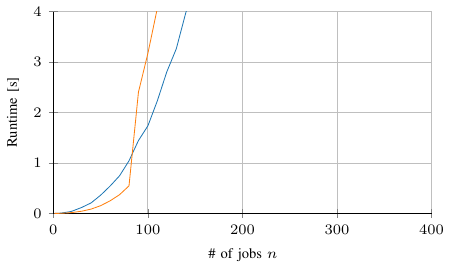}
    \label{fig:0_60_full}
    }
    \subfloat[\centering 15m time granularity]{
    	\includegraphics{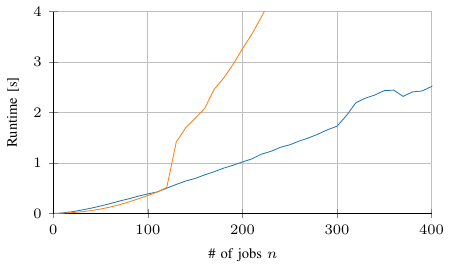}
    \label{fig:0_900_full}
    } \\
    \subfloat[\centering 30m time granularity]{
		\includegraphics{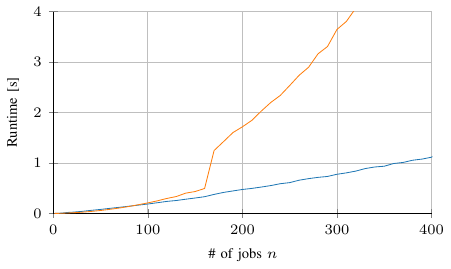}
    \label{fig:0_1800_full}
    }
    \subfloat[\centering 1h time granularity]{
		\includegraphics{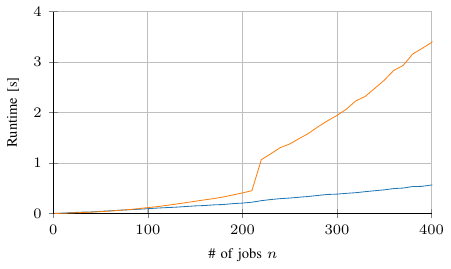}
    \label{fig:0_3600_full}
    }
\caption{Runtimes relative to instance size, for various time granularities, solved for the whole day and using \texttt{shortest\_augmenting\_path()}.}
\label{fig:0_full}
\end{figure*}

In the following, we compare \focs\ to the \qp\ implementation for comparison with (commercial) optimization tools already available for \mpc s. The total runtime (i.e., the sum of the model building and solving), relative to instance size $n$ for the four time granularities is depicted in Figure~\ref{fig:0_full}. Note that for \focs, we assume that the flow method \texttt{shortest\_augmenting\_path()} is applied. All instances are solved for the full day.

First, we note that across subfigures, the quadratic program solved in Gurobi initially either outperforms \focs\ or runs at similar speed. However, we observe a gradually steeper increase in runtime relative to instance size for \qp\ compared to the \focs\ model. As a result, the median runtime of \focs\ outperforms \qp\ starting from instance sizes of approximately 90, 120, 90 and 80 for time granularities of 1~minute, 15~minutes, 30~minutes and 1~hour respectively.
Furthermore, we observe a sudden steep increase in the runtime of the \qp\ implementation, at instance sizes from around 90, 130, 170 and 230 in Figures~\ref{fig:0_60_full}, \ref{fig:0_900_full}, \ref{fig:0_1800_full} and \ref{fig:0_3600_full} respectively. This behaviour possibly reflects that a certain memory threshold is reached after which solving becomes more costly for Gurobi due to swapping memory.

Overall, as expected, a finer time granularity for a given instance size results in a larger runtime for both the \qp\ and \focs. Notably, the slope of the runtime for \focs\ seems to be almost linear for the instance sizes considered, except in Figure~\ref{fig:0_60_full}. There, a non-linear increase can be observed. 
To give an indication of this growth, the median runtimes for \focs\ and \qp\ with instance size 400 under time granularity 1~minute are respectively 36 and 62 seconds. 

\begin{figure*}[t]
\begin{center}
    \includegraphics{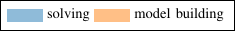}
\end{center}
\centering
    \subfloat[\centering \focs]{
        \includegraphics{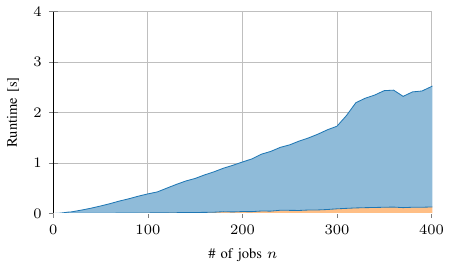}
    \label{fig:0_900_mbvssol_FOCS}
    }
    \subfloat[\centering \qp]{
    	\includegraphics{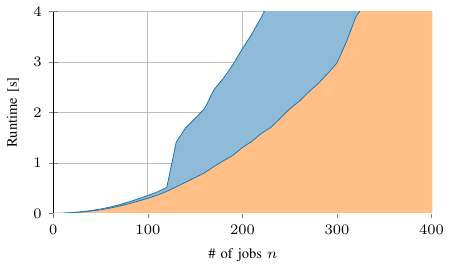}
    \label{fig:0_900_mbvssol_QP}
    }
\caption{Runtimes disaggregated into model building and solving for 15~m granularity.}
\label{fig:0_900_mbvssol}
\end{figure*}

In a next step, we disaggregate the runtimes into model building and solving components (see Figure~\ref{fig:0_900_mbvssol}). 
%
Figure~\ref{fig:0_900_mbvssol_FOCS} presents the results for \focs, whereas Figure~\ref{fig:0_900_mbvssol_QP} presents the results for \qp. Both figures are stackplots, depicting the total runtime as sum of the model building (bottom part, orange shading) and solving (top part, blue shading). The results correspond to experiments with time granularity of 15~minutes and \texttt{shortest\_augmenting\_path()} applied within \focs. 

Notably, \focs\ spends most of its computation time on solving. The almost negligible contribution of model building can be attributed to the mostly dictionary-based implementation of flows in the \texttt{networkx} package. For \qp, on the other hand, model building makes up a considerable portion of the runtime. 
Moreover, Figure~\ref{fig:0_900_mbvssol_QP} clearly attributes the earlier observed sudden steep increase in the runtime of the \qp\ (discussed above) to its solving. This can be seen in the continued smooth increase in the bottom (orange) area, whereas approximately at instance size 130 the contribution of the solving time to the total runtime increases drastically.

\section{Discussion}\label{sec:discussion}
In the previous section, we have presented the empirical results of our experiments. Next, we want to shortly discuss some of the limitations of this research, and put the results in perspective to physical parking lots for EV charging with \mpc. 

First, we consider the implementation of \focs\ as used for the experiments. 
While the overall \focs\ method builds on augmentation of flows, within the used implementation the number of considered intervals is reduced between consecutive iterations and a new maximum flow with the empty flow as initialization is executed.
A more efficient implementation would use the flow from the previous iteration, and augment on top of that. Especially if \focs\ uses a maximum flow method that relies on augmenting paths, this would speed up the optimization. 

\begin{figure}[]
    \centering
    \begin{tabular}{c|c c c c}
         $j$ & $a_j$ & $d_j$ & $e_j$ & $p_j^{\max }$ \\  \hline
         1 & 0 & 2 & 2 & 2 \\
         2 & 0 & 1 & $\frac{1}{2}$ & 2 \\
         3 & 1 & 2 & $\frac{1}{2}$ & 2 \\
         4 & 0 & 2 & 2 & 2 \\
   \end{tabular}

    \includegraphics{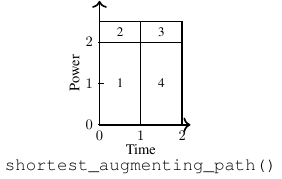}
    \includegraphics{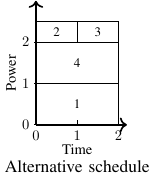}
    \caption{Example of an instance where \texttt{shortest\_augmenting\_path()} results in a schedule with bias for one job over the other.}
    \label{fig:0_unfairSchedule}
\end{figure}

Moreover, in \mpc\ settings, both efficiency and quality of the optimization are crucial.
Therefore, it is important to note that for the conclusions of Section~\ref{sec:resultsMaxFlowComparison} we have not considered the exact form of the found schedule but focused solely on the runtime of the various flow methods.
However, in the setting of EV scheduling, the fairness, robustness and user acceptability of the found schedule itself are of great importance. To illustrate, consider the example in Figure~\ref{fig:0_unfairSchedule}. 
%
Here, jobs~1 and 4 share the same properties. The example provides two schedules. In the schedule on the left those two jobs are scheduled consecutively in disjoint intervals, whereas on the right they are scheduled in parallel. Both schedules are optimal solutions to the problem with constraints (\ref{eq:MIP}) minimizing objective function (\ref{eq:L2NormObjective}). One may argue, however, that the schedule on the right is fairer to the users and more robust to early departure of either job. 
Such considerations were left out of scope of this work, but should be considered in future work considering optimization in \mpc\ settings.

Finally, in the comparison of \focs\ to \qp, we do not take the adaptivity of the models in an \mpc\ setting into account. While the model building for \qp\ takes considerable time, Gurobi is known for efficiently altering and re-solving its models, using previous results in its initialization. This was left out of scope for the comparison. While we do recommend future work in that direction, we first suggest to conduct theoretical research on the efficient initialization of \focs\ and to implement the subsequent results to allow for a reasonable comparison.

\section{Conclusion}\label{sec:conclusion}
This work investigated potential algorithms for an optimization subroutine in model predictive controllers for EV scheduling. To this end, we introduced a pre-mature stop feature for \focs, an algorithm suitable for this subroutine. The proposed feature can decrease runtime in \mpc\ settings. It is particularly effective close to peak hours, for example around noon at an office parking lot when parking lot occupancy is expected to be at its peak. 
Using a 15~minute granularity, the proposed variant of \focs\ shows a median runtime improvement of on average 28\% compared to the original unadapted version to calculate a charging schedule at noon. 
Moreover, this work has illustrated the considerable impact of the chosen maximum flow method for \focs. 
Finally, for large EV parking lots, the presented runtimes show the competitiveness of \focs\ with commercial state-of-the-art optimization software. 

Future work may dive deeper into the integration of \focs\ into an \mpc, and into the properties of the found schedules depending on the used maximum flow method.

\section*{Acknowledgements}
This research is conducted within the SmoothEMS met GridShield project subsidized by the Dutch ministries of EZK and BZK (MOOI32005).

\bibliographystyle{IEEEtran}
\bibliography{literature_references_and_summaries}

\end{document}